\documentclass{amsart}
\theoremstyle{plain}
\newtheorem{theorem}{Theorem}[section]
\newtheorem{lemma}[theorem]{Lemma}
\newtheorem{proposition}[theorem]{Proposition}
\newtheorem{corollary}[theorem]{Corollary}
\theoremstyle{definition}

\newtheorem{example}[theorem]{Example}

\theoremstyle{remark}
\newtheorem{remark}[theorem]{Remark}
\newcommand{\abs}[1]{\vert #1 \vert}
\newcommand{\norm}[1]{\Vert #1 \Vert}
\newcommand{\C}{{\mathbb C}}
\newcommand{\R}{{\mathbb R}}
\newcommand{\N}{{\mathbb N}}

\newcommand{\Z}{{\mathbb Z}}
\newcommand{\SI}{{\rm\Sigma}}

\newcommand{\zp}{\lbrace}
\newcommand{\zz}{\rbrace}
\newcommand{\pmz}{\subseteq}
\newcommand{\D}{{\rm\Delta}}
\newcommand{\bD}{\partial\D}
\newcommand{\bdr}{\partial}
\begin{document}
\title[Riemann-Hilbert problem for Riemann surfaces]
{Nonlinear Riemann-Hilbert problem for bordered Riemann surfaces}
\author{Miran \v{C}erne}
\address{Department of Mathematics\\
University of Ljubljana\\
Jadranska 19, 1\,111 Ljubljana, Slovenia}
\email{miran.cerne@fmf.uni-lj.si}
\thanks{The author was supported in part by a grant from the Ministry of Science
of the Republic of Slovenia.}
\keywords{Riemann-Hilbert problem, bordered Riemann surface, polynomial hull}
\subjclass{Primary 30E25, 35Q15; Secondary 32E20}
\begin{abstract}
Let $\SI$ be a bordered Riemann surface with genus $g$ and $m$ boundary components. Let 
$\zp\gamma_{z}\zz_{z\in\bdr\SI}$ be a smooth family of smooth Jordan curves in $\C$
which all contain the point $0$ in their interior. Then there exists a holomorphic
function $f(z)$ on $\SI$ smooth up to the boundary with at most $2g+m-1$ zeros on 
$\SI$ such that $f(z)\in\gamma_z$ for every $z\in\bdr\SI$.
\end{abstract}
\maketitle
\section{Introduction}
\label{intro}
Let $\SI$ be (the interior of) a bordered Riemann surface with genus $g$ and $m$ real analytic
boundary components. In this paper we consider the existence of holomorphic solutions
of a nonlinear Riemann-Hilbert problem for $\SI$. Let $k$ be a nonnegative integer and let 
$0<\alpha<1$. We denote by $C^{k,\alpha}(\bdr\SI)$ the H\"{o}lder space of all real 
$k$ times differentiable functions on the boundary $\bdr\SI$ whose derivatives of order $k$ are
H\"{o}lder continuous of order $\alpha$ and we denote by 
$A^{k,\alpha}(\SI)$ the space of all holomorphic functions on $\SI$ which
are of class $C^{k,\alpha}$ on $\overline{\SI}$.

We will say that a family of simple closed curves
$\zp\gamma_z\zz_{z\in\SI}$ in $\C$ is a $C^k$ $(k\in\N)$ family of Jordan curves
if there exists a function $\rho\in C^k(\bdr\SI\times\C)$ such that 
$$\gamma_z = \zp w\in\C ; \rho(z,w) = 0\zz$$
and $(\overline{\partial}_w\rho)(z,w)\ne 0$ for every $z\in\bdr\SI$ and $w\in\gamma_z$.
We call $\rho$ a defining function for the family $\zp \gamma_z\zz_{z\in\bdr\SI}$.

\begin{theorem}\label{main}
Let $\SI$ be a bordered Riemann surface with genus $g$ and $m$ real analytic
boundary components. Let $\zp\gamma_{z}\zz_{z\in\bdr\SI}$ be a $C^{k+1}$ $(k\ge 3)$ family 
of Jordan curves in $\C$ which all contain the point $0$ in their interior.
Then there exists a holomorphic function $f\in A^{k,\alpha}(\SI)$
with at most $2g+m-1$ zeros on $\SI$ such that $f(z)\in\gamma_z$ for every $z\in\bdr\SI$.
\end{theorem}

The {\it interior} of a simple closed curve $\gamma\subseteq\C$ is the bounded component of 
$\C\setminus\gamma$. 

\begin{corollary}\label{divisor}
Let ${\mathcal D}\ge 0$ be a divisor of finite degree on $\SI$. Then there exists a
solution $f\in A^{k,\alpha}(\SI)$ of the Riemann-Hilbert problem for $\zp\gamma_{z}\zz_{z\in\bdr\SI}$ 
such that $(f)\ge{\mathcal D}$
and that the degree of the divisor $(f)$ is at most ${\rm deg}({\mathcal D})+ 2g + m - 1$.
\end{corollary}

In the disc case $(\SI=\D)$ results of this type were first proved by \v{S}nirelman, \cite{sni},
and later by Forstneri\v{c}, \cite{for1}, in a different context. For planar domains 
the result, without the bound on the number of the zeros,
was proved by Begehr and Efendiev in \cite{beg-efe}. Linear Riemann-Hilbert problem for bordered Riemann
surfaces was considered by Koppelman in \cite{kop}. For other results on the 
existence of solutions for planar domains one should see 
\cite{efe-wen1,efe-wen2,efe-wen3,mak-efe,weg} and the references therein.

Let 
$$T=\bigcup_{z\in\bdr\SI} (\zp z\zz\times\gamma_z) .$$
Then $T$ is the union of $m$ totally real tori in $\bdr\SI\times\C$ and the existence
of a solution of the Riemann-Hilbert problem for $\zp \gamma_z\zz_{z\in\bdr\SI}$
is equivalent to the existence of a holomorphic
function $f\in A^{k,\alpha}(\SI)$ whose graph is attached to $T$,
that is, $(z,f(z))\in T$ for every $z\in\bdr\SI$. Results
on the existence of large analytic discs attached to certain manifolds can be found
in \cite{for1,gro,slo3,whi2,whi3,whi4}. Results on small holomorphic perturbations
of a given analytic object along maximal real manifolds are given in
\cite{cer1,cer2,cer-glo,cer-for,for2,glo1,glo2}. 

There is also a strong connection of the solutions of the Riemann-Hilbert problem 
to the description of the polynomial hull of a compact set $K$ fibered over $\bdr\D$, 
\cite{ale-wer1,cer3,cer4,for1,law1,law2,slo1,slo2,slo3,whi1,whi2,whi3,whi4}.
Recall that the polynomial hull $\widehat{K}$ 
of a compact set $K\pmz\C^n$ is defined as 
$$\widehat{K} = \zp z\in\C^n\, ; \abs{p(z)}\le\max_K\abs{p}
{\rm\ for\ every\ polynomial\ } p {\rm\ in\ } n {\rm\ variables}\zz$$
and that by the maximum principle every analytic variety whose boundary 
lies in $K$ belongs to the polynomial hull $\widehat{K}$ of $K$.

\begin{corollary}\label{hull}
Let $a_0,\dots,a_{n-1}$ be holomorphic functions on the disc $\D$ which are smooth up to 
the boundary. Let
$$V = \zp (z,w)\in\overline{\D}\times\C; w^n + a_{n-1}(z) w^{n-1} + \dots + a_0(z) = 0\zz$$
be an $n$-sheeted analytic variety over $\overline{\D}$ such that for each $z\in\bdr\D$ there are exactly
$n$ different points in $\zp z\zz\times\C$ which belong to $\partial V$. Let $\zp \gamma_{p}\zz_{p\in\bdr V}$ 
be a $C^{k+1}$ $(k\ge 3)$ family of Jordan curves in $\C$ such that each $\gamma_p$ contains the point $w(p)$, 
$p=(z(p),w(p))$, in its interior. Then there exists an $n$-sheeted analytic variety $V_0$ over $\D$ 
whose boundary is contained in 
$$\bigcup_{p=(z(p),w(p))\in\bdr V} (\zp z(p)\zz\times\gamma_p) .$$
\end{corollary}

The main theorem has an application to embedding finitely
connected planar domains into convex domains in $\C^n$ $(n\ge 2)$. See
\cite{cer-glo,cer-for} for embeddings into $\C^2$. 
\begin{corollary}\label{embedding}
Let $\Omega$ be a smoothly bounded convex domain in $\C^n$ $(n\ge 2)$ and let $D$
be a bounded finitely connected planar domain with smooth boundary. Then there exists a holomorphic
embedding of $D$ into $\Omega$.
\end{corollary}

The paper is organized as follows. In Section \ref{RH-problem} the proofs of the main theorem and corollaries
are given. In Section \ref{sec-example} a special case where $\zp\gamma_z\zz_{z\in\bdr\SI}$ 
are circles centered at $0$ and $\SI$ is a planar domain is considered more closely.
Finally in the Appendix some technical results needed in our arguments are proved.

\section{The Riemann-Hilbert problem}
\label{RH-problem}

Our first proposition claims that if there is a solution of the Riemann-Hilbert problem, then
there is also one which has at most $2g+m-1$ zeros on $\SI$. The idea of the proof is the fact that
the family of solutions of the Riemann-Hilbert problem with more than $2g+m-1$ zeros is not compact
in the space $A^{k,\alpha}(\SI)$. 

\begin{proposition}\label{zeros}
If there exists a solution $f\in A^{k,\alpha}(\SI)$ of the Riemann-Hilbert problem for $\zp \gamma_z\zz_{z\in\bdr\SI}$, 
then there exists a solution with at most $2g+m-1$ zeros.
\end{proposition}
\begin{proof}
Every solution $f\in A^{k,\alpha}(\SI)$ of the Riemann-Hilbert problem for $\zp \gamma_z\zz_{z\in\bdr\SI}$
has finitely many zeros on $\SI$. 
Let $N$ be the minimal number of zeros counting algebraically of a solution of the Riemann-Hilbert problem 
and let $f_0\in A^{k,\alpha}(\SI)$ be a solution with the minimal number of zeros. 
Let ${\mathcal F}$ be the family of all solutions $f\in A^{k,\alpha}(\SI)$ of the Riemann-Hilbert
problem such that $f : \bdr\SI\rightarrow\C\setminus\zp 0\zz$ is homotopic to 
$f_0 : \bdr\SI\rightarrow\C\setminus\zp 0\zz$, that is, $f$ and $f_0$ have the same winding number
over each boundary component of $\SI$ (we orient the boundary $\partial\SI$ coherently with the natural
orientation on $\SI$). The number of zeros of a function $f\in A^{k,\alpha}(\SI)$
which is nonzero on $\bdr\SI$ equals the winding number of $f$ along $\bdr\SI$.
Hence each $f\in{\mathcal F}$ has the same number of zeros counting algebraically.

Let $T = \cup_{z\in\bdr\SI} (\zp z\zz\times\gamma_z)$. Manifold $T$ is a finite union of maximal real tori
in $\partial\SI\times\C$. By Lemma \ref{symplectic} in the Appendix there exists a $C^{k+1}$ strongly 
plurisubharmonic function $v$ on $\overline{\SI}\times\C$ such that $T$ is Lagrangian for the symplectic form 
$\omega = i \partial\overline{\partial} v$ and that the $\omega$-area of any disc 
$\zp z\zz\times\widehat{\gamma_z}$ $(z\in\bdr\SI)$ is $1$. Then, by the Stokes' theorem, all the
graphs of the functions in ${\mathcal F}$ have $\omega$-area equal to $N+a$, where $a$ is a constant
which depends only on $T$ and $\omega$.

Using Gromov's compactness theorem (see \cite{gro}, \cite{aud-laf} or \cite{ye}), 
the family ${\mathcal F}$ is compact.
Namely, if $\zp f_n\zz\subseteq{\mathcal F}$ is a sequence, then there exists a subsequence $\zp f_{n_j}\zz$,
a finite set $\Gamma\subset\bdr\SI$ and a holomorphic function $f_{\infty}\in A^{k,\alpha}(\SI)$ such that
$f_{\infty}$ is a solution of the Riemann-Hilbert problem and the sequence $\zp f_{n_j}\zz$ converges to $f_{\infty}$ in
the $C^{k,\alpha}$ sense on the compact subsets of $\overline{\SI}\setminus\Gamma$. Also, at the points $p\in\Gamma$ 
holomorphic bubbles appear
$$\varphi_p : (\D,\bD)\rightarrow ({\rm Int}(\widehat{\gamma_p}), \gamma_p) \ \ \ \ (p\in\Gamma)$$
and the sum of the $\omega$-area of the graph of $f_{\infty}$ and the $\omega$-areas of the bubbles $\varphi_z$,
$z\in\Gamma$, equals the $\omega$-areas of the graphs of $f_{n_j}$'s.
Hence 
$${\rm the\ number\ of\ zeros\ of\ } f_{\infty} + {\rm card\ }\Gamma = 
{\rm the\ number\ of\ zeros\ of\ } f_{n_j} = N .$$
Since $N$ is the minimal number of zeros of a solution of the Riemann-Hilbert problem on $\SI$, we
conclude that $\Gamma$ is empty and therefore ${\mathcal F}$ is a compact subset of $A^{k,\alpha}(\SI)$.

Let $d$ be a distance-metric on $\overline{\SI}$ induced from a Riemannian metric on $\overline{\SI}$.
For each $f\in{\mathcal F}$ let $Z_f\subseteq\SI$ denote the set of its zeros. We define
$$\sigma(f) = d(Z_f, \bdr\SI) \ \ ({\rm the\ distance\ of\ } Z_f {\rm \ from\ } \bdr\SI)\ \ \ \
{\rm\ and}\ \ \ \ \sigma = \inf_{f\in{\mathcal F}} \sigma(f) .$$
Then there exists a sequence $\zp f_n\zz\subseteq{\mathcal F}$ such that $\lim_{n\rightarrow\infty} \sigma(f_n) = \sigma$. 
Let $\zp f_{n_j}\zz$ be a subsequence which converges in the $C^{k,\alpha}$ sense on $\overline{\SI}$ to
a function $f_{\infty}\in{\mathcal F}$ (compactness!). 
Since $\sigma(f_{\infty})>0$ and since every function from ${\mathcal F}$ has the same number of 
zeros counting algebraically we get $\sigma(f_{\infty}) = \sigma > 0$. 

Let $\rho$ be a $C^{k+1}$ defining function for 
$\zp\gamma_z\zz_{z\in\bdr\SI}$. The map
$$\Psi: A^{1,\alpha}(\SI)\longrightarrow C^{1,\alpha}(\bdr\SI)$$
$$(\Psi (f))(z) = \rho (z, f(z))$$
is (twice) differentiable, \cite{hil-tai},
and its derivative at a function $f\in A^{1,\alpha}(\SI)$ is
$$(D\Psi(f)h)(z) = 2{\rm Re}((\partial_w\rho)(z,f(z))h(z)) .$$
Let $f$ be a solution of the Riemann-Hilbert problem and 
let $\kappa$ be the winding number of the nonzero function
$(\overline{\partial}_w\rho)(z,f(z))$ on $\bdr\SI$. Geometric assumptions on
the family $\zp \gamma_z\zz_{z\in\bdr\SI}$ imply that $\kappa$ equals
the winding number of $f$ on $\bdr\SI$ and hence the number of its zeros.
The linear operator $D\Psi(f)$ is a Fredholm operator of index 
$2\kappa - (2g+m-2)$ and it has no cokernel if $\kappa\ge 2g + m - 1$, \cite{kop}.
Therefore the implicit function theorem in Banach spaces shows
(see also \cite{cer-glo,cer-for}) that in the case $N\ge 2g + m - 1$
the family ${\mathcal F}$ is $M = 2N - (2g+m-2)$ dimensional submanifold of $A^{1,\alpha}(\SI)$. 
Using a theorem of \v{C}irka, \cite{chi}, we may even conclude that ${\mathcal F}$ is
an $M$ dimensional submanifold of $A^{k,\alpha}(\SI)$.

Let us now assume that $N\ge 2g+m$. 
Let $f_0\in{\mathcal F}$ and $z_0\in\SI$ be such that $f_0(z_0)=0$ and $\sigma = d(z_0,\bdr\SI)$. 
Let $f(\cdot,p)$ be a local $C^1$ parametrization of ${\mathcal F}$ with the space 
${\mathcal P}= \R^M$ such that $f(\cdot,0) = f_0$. The derivative $(D_p f)(\cdot,0)$ is an isomorphism
from ${\mathcal P}$ into the tanget space of ${\mathcal F}$ at the point $f_0$, that is, for every
$p\in{\mathcal P}$ we have
$${\rm Re}(\overline{a(z)}\, (D_p f)(z,0)p)=0\ \ \ \ (z\in\bdr\SI) ,$$
where $a(z) = (\overline{\partial}_w\rho)(z,f_0(z))$.

We define the mapping 
$\Phi : \SI\times{\mathcal P}\rightarrow\SI\times\C$ by $\Phi(z,p)= (z, f(z,p))$.
We will show that $\Phi$ is a submersion in a neighbourhood of $(z_0,0)\in\SI\times{\mathcal P}$. Hence for 
each point $z$ close enough to $z_0$ the equation $\Phi(z,p) = (z,0)$ has a solution $p(z)$
and this contradicts the fact that $z_0$ is a zero of a function from ${\mathcal F}$ 
which is the nearest to the boundary $\bdr\SI$. Therefore $N\le 2g+m-1$. 

To prove that $\Phi$ is a submersion in a neighbourhood of $(z_0,0)$ we have to prove that
the derivative $(D\Phi)(z_0,0) : T_{z_0}\SI\times{\mathcal P}\rightarrow T_{z_0}\SI\times\C$ is surjective.
For this it is enough to prove that the partial derivative $(D_p f)(z_0,0) : {\mathcal P}\rightarrow\C$ is surjective. 

Let us assume that the image of $(D_p f)(z_0,0)$ is either $0$ or $1$ dimensional. In either case 
its image lies in a line in $\R^2$, that is, there exists a nonzero complex number $a_0$ such that 
$${\rm Re}(\overline{a_0}\, (D_p f)(z_0,0)p) = 0$$
for every $p\in{\mathcal P}$. We may assume, without loss of generality, that $a_0 = i$.

Let $h_j(z) = (D_p f)(z,0) e_j$, where $e_j$
$(j=1,\dots,M)$ is a basis of the space ${\mathcal P}$.
Then $h_j(z_0)\in\R$ for each $j$ and hence there is another basis $\widetilde{e}_j$
$(j=1,\dots, M)$ of the space ${\mathcal P}$ such that for 
$\widetilde{h}_j(z) = (D_p f)(z,0)\widetilde{e}_j$ we have $\widetilde{h}_j(z_0) = 0$
for every $j=1,\dots, M-1$.

We know that every $(D_p f)(z,0) p$  is a solution of the homogeneous Riemann-Hilbert problem
$${\rm Re}(\overline{a(z)}\, (D_p f)(z,0)p)=0\ \ \ \ (z\in\bdr\SI).$$
Hence
$${\rm Re}(\overline{a(z)}\, \widetilde{h}_j(z))=0$$
for every $z\in\bdr\SI$ and $j=1,\dots, M$. Let $g$ be a function from $A^{k,\alpha}(\SI)$ such that
it has the only simple zero on $\overline{\SI}$ at the point $z_0$. Thus
$\widetilde{h}_j(z) = g(z)\, g_j(z)$ for some $g_j\in A^{k,\alpha}(\SI)$ and $j=1,\dots, M-1$.
Also, since $\widetilde{h}_j$, $j=1,\dots,M-1$, are linearly independent, so are the functions $g_j$,
$j=1,\dots,M-1$.
Therefore
$${\rm Re}(\overline{a(z)}\, g(z)\, g_j(z))=0$$
for every $z\in\bdr\SI$ and $j=1,\dots,M-1$. Let us define
$$ b(z) = a(z)\, \overline{g(z)}$$
for $z\in\bdr\SI$. Then $b$ is a nonzero function on $\bdr\SI$ of class $C^{k}$ and its winding
number on $\bdr\SI$ is
$$W(b) = W(a) - W(g) = W(a) - 1 = N - 1.$$
Hence $W(b)\ge 2g + m - 1$ and the space of solutions of the linear homogeneous Riemann-Hilbert
problem ${\rm Re}(\overline{b(z)}\, h(z))=0$ is $2(N-1) - (2g+m-2) = M - 2$ dimensional, \cite{kop}. Therefore
the functions $g_j$, $j=1,\dots, M-1$ are linearly dependent, which is a contradiction.	
\end{proof}

\begin{remark}
The above proof shows that in the case the minimal number of zeros of a solution of the Riemann-Hilbert
problem on $\SI$ is $2g+m-1$, then the set of solutions with the minimal number of zeros 
is a closed submanifold in the space $A^{k,\alpha}(\SI)$ of dimension $2g+m$. See also Example \ref{example}.
\end{remark}

The idea of the proof of the main theorem is the observation by Begehr and Efendiev, \cite{beg-efe},
that one can use solutions of the Riemann-Hilbert problem on the disc to get a very good approximate
solutions of the original Riemann-Hilbert problem. Namely, for each boundary component $\Gamma$
of $\SI$ one can get a holomorphic function on its neighbourhood in $\SI$ (an annulus)
such that this function solves the Riemann-Hilbert problem on $\Gamma$ and it is uniformly
close to $0$ away from $\Gamma$. All these functions are then extended to $\SI$ so that
their $\overline{\partial}$ derivatives are very small. Using solutions of the $\overline{\partial}$
equation on $\SI$ we get a holomorphic function on $\SI$ whose boundary values almost solve
the original Riemann-Hilbert problem. Finally we use the following theorem from
\cite{beg-efe} to complete the proof.
\begin{theorem}[Begehr-Efendiev] \label{Begehr-Efendiev}
Let $X$ and $Y$ be Banach spaces. Let $A : X\rightarrow Y$ be a continuous mapping,
Fr\'{e}chet differentiable in a neighbourhood of $x_0\in X$. Assume that $(DA)(x_0)$
has a bounded right inverse $B(x_0)$ with the norm $\omega_1>0$. Also,
let $\omega_2>0$ be such that the derivative $(DA)(x)$ satisfies a Lipschitz condition
$$\norm{(DA)(x_1)-(DA)(x_2)}<\omega_2 \norm{x_1-x_2}$$ 
for $x_1, x_2$ in a neighbourhood of $x_0$. Finally, let $\omega_3>0$ be such that 
$$4 \omega_1 (\omega_1 + 1)(\omega_2 + 1)\omega_3<1.$$
Then for every $y\in Y$ such that $\norm{A(x_0)-y}<\omega_3$ there exists a solution of the equation 
$A(x) = y$. 
\end{theorem} 

\begin{proof}(Theorem \ref{main}) Let $\rho\in C^{k+1}(\bdr\SI\times\C)$ be a defining function
for $\zp\gamma_z\zz_{z\in\bdr\SI}$ and let $\Gamma_1$,$\dots$,$\Gamma_m$ denote the boundary components 
of $\SI$. Let $A_1,\dots,A_m$ be pairwise disjoint annuli in $\SI$ with real analytic boundaries
such that $\Gamma_j\subset\bdr A_j$ is one boundary component of $A_j$ $(j=1,\dots,m)$. 
For $j=1,\dots, m$ there exists $0<r_j<1$ and an up to the boundary real analytic
biholomorphism
$$\Phi_j: (\D(0,1)\setminus\overline{\D(0,r_j)}, \bD)\longrightarrow (A_j,\Gamma_j).$$
Now we can solve the Riemann-Hilbert problem for $\D(0,1)$ and 
$\zp \gamma_{\Phi_j(z)}\zz_{z\in\bD(0,1)}$ as the boundary data, \cite{for1}. 
Even more, the proof in \cite{for1} (see Lemma \ref{estimates} in the Appendix) 
implies that there exists an integer $l$ such that for any $n\in\N$ and any compact subset 
$\overline{\D(0,q)}\subseteq\D(0,1)$ $(0<q<1)$
there exists a solution $\widetilde{f}_{nj}$ of the Riemann-Hilbert problem on $\D$
with the winding number $W(\widetilde{f}_{nj})=n$ and such that
$$ \norm{\widetilde{f}_{nj}}_{A^{1,\alpha}(\D(0,1))}\le C n^l\ \ \ \ {\rm and}\ \ \ \ 
\norm{\widetilde{f}_{nj}}_{A^{1,\alpha}(\D(0,q))}\le C n^l q^{n-2}$$
for some finite positive constant $C$ which does not depend on $n$ and $j$.

Let $\max\zp r_1,\dots,r_m\zz<q<1$ and let $K_j = \Phi_j(\overline{\D(0,q)}\setminus \D(0,r_j))$
be a compact subset of $\overline{A_j}\setminus\Gamma_j$, $(j=1,\dots,m)$.
Then there exists a sequence of holomorphic functions $f_{nj} = \widetilde{f_{nj}}\circ\Phi_j^{-1}$ on $A_j$
such that $f_{nj}(z)\in\gamma_z$ for each $z\in\Gamma_j$, 
the winding number of $f_{nj}$ along $\Gamma_j$ is $n$, and 
\begin{equation}\label{piecewise}
\norm{f_{nj}}_{A^{1,\alpha}(A_j)}\le C n^l\ \ \ \ {\rm and}\ \ \ \norm{f_{nj}}_{C^{1,\alpha}(K_j)}\le C n^l q^{n}
\end{equation}
for some constant $C<\infty$ which does not depend on $n$ and $j$. 

Let $\chi$ be a smooth function on $\overline{\SI}$, $0\le \chi\le 1$,  such that $\chi = 1$ on
a neighbourhood of $\cup_{j=1}^m(\overline{A_j\setminus K_j})$ and such that the support of $\chi$ is contained in 
$\cup_{j=1}^m (A_j\cup\Gamma_j)$.
For each $n\in\N$ we define $f_n : \overline{\SI}\rightarrow\C$ as
$$ f_n = \sum_{j=1}^m \chi\,f_{nj} .$$
Then 
$$\norm{f_n}_{C^{1,\alpha}(\SI)}\le C\,n^l\ \ \ \ {\rm and}\ \ \ \ \norm {\overline{\partial} f_n}_{C^{\alpha}(\SI)}\le C n^l q^{n} .$$
Therefore there exists a solution $u_n$ on $\SI$ of the equation $\overline{\partial} u_n = \overline{\partial} f_n$
such that $\norm {u_n}_{C^{1,\alpha}(\SI)}\le C n^l q^{n}$ for some universal constant $C$, \cite{kop,vek}. 
Hence for the functions $h_n = f_n - u_n\in A^{1,\alpha}(\SI)$ we have
$$\norm{\rho(z,h_n(z))}_{C^{1,\alpha}(\bdr\SI)}\le C n^{3l} q^{n}.$$ 
See Lemma \ref{approximation} in the Appendix.

We want to use Theorem \ref{Begehr-Efendiev}. Let us observe the twice differentiable map 
$\Psi: A^{1,\alpha}(\SI)\rightarrow C^{1,\alpha}(\bdr\SI)$ defined as $(\Psi(f))(z) = \rho(z,f(z))$.
The arguments above show that the initial error for $\Psi$ at the point $h_n$ is $\omega_3 = C n^{3l} q^n$.   
Now we have to find an estimate of the norm of a right inverse of the 
derivative of $\Psi$ at the points $h_n$ and an estimate on the derivative's Lipschitz constant near $h_n$. 

The derivative of $\Psi$ at $h_n$ is of the form
$$(\Psi_n(h))(z) = 2{\rm Re}((\partial_w\rho)(z,h_n(z))h(z)) .$$
We define
$$(\widetilde{\Psi_n}(h))(z) = 2{\rm Re}((\partial_w\rho)(z,f_n(z))h(z))$$
a linear mapping from $A^{1,\alpha}(\SI)$ into $C^{1,\alpha}(\bdr\SI)$.
Then from Lemma \ref{compositum} we get
$$\norm{\widetilde{\Psi_n}}\le 2\norm{(\partial_w\rho)(z,f_n(z))}_{1,\alpha}
\le C\,(\norm{f_n}_{1,\alpha}^2+1)\le C n^{2l}$$
and
$$\norm{\Psi_n - \widetilde{\Psi_n}} \le 2\norm{(\partial_w\rho)(z,h_n(z)) - (\partial_w\rho)(z,f_n(z))}_{1,\alpha}\le$$
$$\le C\, (\norm{f_n}_{1,\alpha}^2+\norm{h_n-f_n}_{1,\alpha}^2+1)\, \norm{h_n - f_n}_{1,\alpha}\le C n^{3l} q^n$$
as in the proof of Lemma \ref{approximation} in the Appendix.
Hence, if we can find a right inverse $\widetilde{B_n}$ of $\widetilde{\Psi_n}$
such that $\norm{\widetilde{B_n}}\le C n^{l_1}$ for some fixed $l_1\in\N$ and $C>0$, then 
$$B_n = \widetilde{B_n}\, (I + (\Psi_n - \widetilde{\Psi_n})\widetilde{B_n})^{-1}$$
is a right inverse of $\Psi_n$ such that $\norm{B_n}\le C n^{l_1}$ for $n$
large enough.

Let $g\in C^{1,\alpha}(\bdr\SI)$. Using solutions of the linear Riemann-Hilbert problem on the disc, 
Lemma \ref{nonhomogeneous}, we find, similarly as in (\ref{piecewise}), that
there exist $l_1\in\N$ and a finite constant $C$ such that
for each boundary component $\Gamma_j\subseteq\bdr\SI$, $j=1,\dots, m$,
there is a holomorphic function $H_{nj}(g)$ on $A_j$ with the properties
$$\norm{H_{nj}(g)}_{A^{1,\alpha}(A_j)} \le C\ n^{l_1}\, \norm{g}_{1,\alpha},\ \ \ \ 
\norm{H_{nj}(g)}_{A^{1,\alpha}(K_j)} \le C\, n^{l_1}\, q^n\, \norm{g}_{1,\alpha}$$
and
$2{\rm Re}((\partial_w\rho)(z,f_n(z))\,H_{nj}(g)) = g$ on $\Gamma_j$.
The construction of the functions $H_{nj}(g)$ is linear in $g$, Lemma \ref{nonhomogeneous}.

We define $H_n(g) = \sum_{j=1}^m \chi\,H_{nj}(g)$ and observe that 
$$\norm{H_n(g)}_{1,\alpha} \le C\, n^{l_1}\, \norm{g}_{1,\alpha}\ \ \ \ {\rm and}\ \ \ \
\norm{\overline{\partial}H_n(g)}_{\alpha}\le C\, n^{l_1}\, q^n\, \norm{g}_{1,\alpha} .$$

Let $T_0 : C^{\alpha}(\overline{\SI})\rightarrow C^{1,\alpha}(\overline{\SI})$ be an operator which gives a 
solution of the $\overline{\partial}$-equation $\overline{\partial} T_0 u = u$
such that $\norm{T_0 u}_{1,\alpha} \le C \norm{u}_{\alpha}$ for some positive constant $C$
and every $u\in C^{\alpha}(\overline{\SI})$, \cite{kop,vek}. Let us define
$$\widehat{B_n}g = H_n(g) - T_0 (\overline{\partial}H_n(g)) .$$
Then $\widehat{B_n}$ is a linear operator from $C^{1,\alpha}(\bdr\SI)$ into $A^{1,\alpha}(\SI)$
with the norm
$$\norm{\widehat{B_n}} \le C\, n^{l_1}$$
and
$$\widetilde{\Psi_n}\widehat{B_n} g = g - \widetilde{\Psi_n} T_0(\overline{\partial}H_n(g)) =
g - E_n(g)$$
where $\norm{E_n}\le C\, n^{2l+l_1}\, q^n$ (the operator $\widetilde{\Psi_n}$ can be naturally
extended to the space $C^{1,\alpha}(\overline{\SI})$ with the same norm estimate).
Hence $\widetilde{B_n}=\widehat{B_n}(I-E_n)^{-1}$ is a left inverse of $\widetilde{\Psi_n}$ such that its norm
is bounded by $C n^{l_1}$ for some universal constant $C$.

At the end we estimate the Lipschitz constant for the derivative $D\Psi$ of the mapping $\Psi$
near the point $h_n$. Here we need that $\rho$ is at least of class $C^4$.
Let $f_0$ and $f_1$ be two functions such that
$\norm{f_0 -h_n}_{1,\alpha} < 1$ and $\norm{f_1-h_n}_{1,\alpha}<1$. 
Then, as in the proof of Lemma \ref{approximation}, we get
$$\norm{D\Psi(f_1)-D\Psi(f_0)}\le 2\norm{(\partial_w\rho)(z,f_1(z))-(\partial_w\rho)(z,f_0(z))}_{1,\alpha}\le$$
$$\le 2\norm{f_1-f_0}_{1,\alpha}\, 
\sup_{t\in\lbrack 0,1\rbrack} \norm{(D^2\rho)(z,t f_1 + (1-t) f_2)}_{1,\alpha}\le$$
$$\le C (\norm{h_n}^2_{1,\alpha} + 1)\, \norm{f_1-f_0}_{1,\alpha}\le C n^{2l}\norm{f_1-f_0}_{1,\alpha} .$$

Hence we have $\omega_1 = C n^{l_1}$, $\omega_2 = C n^{2l}$, $\omega_3 = C n^{3l} q^n$ and for $n$ large enough 
the product $4\omega_1(\omega_1+1)(\omega_2+1)\omega_3$ is arbitrarily small. Therefore by Theorem \ref{Begehr-Efendiev} 
there exists a solution $f\in A^{1,\alpha}(\SI)$ of the Riemann-Hilbert problem
and by \cite{chi} we have $f\in A^{k,\alpha}(\SI)$.
\end{proof}

An immediate corollary of Theorem \ref{main} is

\begin{corollary}
Let $F$ be a continuous function on $\overline{\SI}$ which is holomorphic on $\SI$.
Let $\zp\gamma_{z}\zz_{z\in\bdr\SI}$ be a $C^{k+1}$ $(k\ge 3)$ family 
of Jordan curves in $\C$ such that for each $z\in\bdr\SI$ the point $F(z)$ lies in the interior of $\gamma_z$.
Then there exists a holomorphic function $f\in A^{k,\alpha}(\SI)$
such that $f(z)\in\gamma_z$ for every $z\in\bdr\SI$.
\end{corollary}

Now we give the proofs of corollaries stated in introduction.

\begin{proof}(Corrollary \ref{divisor})
Let $g\in A^{k+1}(\SI)$ be a holomorphic function such that $(g)={\mathcal D}$
and $g(z)\ne 0$ for every $z\in\bdr\SI$.
We define a new $C^{k+1}$ family of Jordan curves in $\C$ enclosing the point $0$ by
$$\widetilde{\gamma}_z = \frac{1}{g(z)}\gamma_z\ \ \ \ (z\in\bdr\SI) .$$
Then there exists a function $\widetilde{f}\in A^{k,\alpha}(\SI)$ such that
$\widetilde{f}(z)\in\widetilde{\gamma}_z$ for every $z\in\bdr\SI$ and such that 
$\widetilde{f}$ has at most $2g+m-1$ zeros on $\SI$. Hence 
$f(z) = g(z) \widetilde{f}(z)$ is the required function. 
\end{proof}

\begin{proof} (Corollary \ref{hull})
Without loss of generality we may assume that $a_0,\dots,a_{n-1}$ are polynomials and that
$V$ is an embedded bordered Riemann surface $\overline{\SI}$ in $\overline{\D}\times\C$. 
If not, we can uniformly on $\overline{\D}$ approximate functions $a_j$ $(j=0,\dots,n-1)$
with polynomials (which we still denote by $a_j$) and replace $V$ by the solution of the equation 
$$w^n + a_{n-1}(z) w^{n-1} + \dots + a_0(z) = a,$$
where $a$ is in absolute value small regular value of the function 
$P(z,w)=w^n + a_{n-1}(z) w^{n-1} + \dots + a_0(z)$.
Let 
$$\Phi = (f,g) : \SI\longrightarrow V\cap (\D\times\C)$$
be a biholomorphism which is smooth up to the boundary. 
For each $p\in\bdr\SI$ let 
$$\widetilde{\gamma}_p = \gamma_{\Phi(p)} - g(p) .$$
Then $\zp \widetilde{\gamma}_p\zz_{p\in\bdr\SI}$ is a $C^{k+1}$ family of Jordan curves
in $\C$ such that $\widetilde{\gamma}_p$ contains the point $0$ in its interior. 
Hence there is a holomorphic function $h\in A^{k,\alpha}(\SI)$ such that 
$h(p)\in \widetilde{\gamma}_p$ for every $p\in\bdr\SI$. We define
$V_0$ as the image of the proper mapping $p\mapsto (f(p),g(p)+h(p))$
from $\overline{\SI}$ into $\overline{\D}\times\C$.	
\end{proof}

\begin{proof}(Corollary \ref{embedding})
Without loss of generality we may assume $n=2$ and that $\overline{D}\times \zp 0\zz\subseteq \Omega$.
Then $\zp z\zz\times\C$ intersects $\bdr\Omega$ transversally for every $z\in\bdr D$ and 
hence $\gamma_{z} = (\zp z\zz\times\C) \cap \bdr \Omega$, $z\in\bdr D$, is a smooth family of smooth Jordan curves in $\C$ which all
contain $0$ in their interior. By Theorem \ref{main} there exists a holomorphic function $f : D\rightarrow\C$
smooth up to the boundary
such that $f(z)\in\gamma_{z}$ for each $z\in\bdr D$. The mapping
$z\mapsto (z,f(z))$ gives a holomorphic embedding of $D$ into $\Omega$. 
\end{proof}

\section{Example}
\label{sec-example}

The starting point for our next observation was the simplest nontrivial possible case one
can try to solve directly by computations: $\gamma_z$ are circles centered at $0$.
Our computations lead to the following result which is a kind of an inversion
problem result for planar domains, e.g. \cite[p.281]{spr}.

Let $D$ be a bounded finitely connected planar domain with smooth boundary. Let 
$\Gamma_0,\dots,\Gamma_m$ be its boundary components, where $\Gamma_0$ is the boundary
of the unbounded component of $\C\setminus D$. Let $h_1,\dots,h_m$ be
harmonic functions on $D$ smooth up to the boundary such that $h_j = 1$ on $\Gamma_j$
and $h_j=0$ on other boundary components. We denote by 
$\overline{D}^m$ the product of $m$ copies of $\overline{D}$ and by $T^m$ the $m$ torus
$S^1\times\dots\times S^1$ ($m$ factors), where $S^1 = \R/\Z$.
\begin{proposition}
The maping $\Phi : \overline{D}^m\longrightarrow T^m$ defined by
$$\Phi(z_1,\dots,z_m) = (\sum_{j=1}^m h_1(z_j),\dots, \sum_{j=1}^m h_m(z_j))\ {\rm mod}\ \Z $$
is surjective.
\end{proposition}
\begin{proof}
Without loss of generality we may assume that 
$$D = \D\setminus\bigcup_{j=1}^m \D(a_j,r_j)$$
where $\abs{a_j}+r_j<1$ for every $j=1,\dots, m$ and $\abs{a_j-a_t}>r_j + r_t$ for every $j\ne t$. 
Let $R(z)$ be any $C^{k+1}$ $(k\ge 3)$ positive function on $\bdr D$. We know that there is $f\in A^{k,\alpha}(D)$
such that $\abs{f(z)} = R(z)$ for every $z\in\bdr D$ and such that it has at most $m$ zeros on $D$. 
Let $z_1,\dots,z_M$ $(M\le m)$ be the zeros of $f$ and let $k_j$ $(j=1,\dots,m)$ be the winding number of
$f$ along $\Gamma_j=\bdr\D(a_j,r_j)$ oriented coherently to $\D(a_j,r_j)$.
Then $f$ is of the form
$$f(z) = (z-z_1)\dots (z-z_M) (z-a_1)^{k_1}\dots (z-a_m)^{k_m} e^{F(z)}$$
for some $F\in A^{k,\alpha}(D)$. 

Let $u_R$ be the harmonic function on $D$ such that $u_R = \log R$ on $\bdr D$
and let $u_j$ $(j=1,\dots,M)$ be the harmonic function on $D$ such that $u_j(z) = -\log\abs{z-z_j}$
on $\bdr D$. Hence
$${\rm Re}(F(z)) = u_R(z) + \sum_{j=1}^M u_j(z) - \sum_{j=1}^m k_j\log\abs{z-a_j}$$
on $\overline{D}$. 

Since $F$ is holomorphic on $D$, for $u={\rm Re}(F)$ we have
\begin{equation}\label{conjugate}
\int_{\Gamma_p}\frac{\partial u}{\partial {\mathbf n}} ds = 0
\end{equation}
for every $p=1,\dots, m$. Here, $\Gamma_p = \bdr\D(a_p,r_p)$ is oriented 
coherently to $\D(a_p,r_p)$ and ${\mathbf n}$ is the outer unit normal to $\bD(a_p,r_p)$.
We know that
$$\int_{\Gamma_p}\frac{\partial \log\abs{z-a_j}}{\partial {\mathbf n}} ds = 
\left\{ \begin{array}{ccl}
2\pi  &;& p=j\\
0     &;& p\ne j \end{array}\right. .$$
Let $G(z,w)$ be the Green's function for $D$. Then we observe that
$$G(z,z_j) = \frac{1}{2\pi}(u_j(z) + \log\abs{z-z_j})$$
and hence
$$\int_{\Gamma_p}\frac{\partial u_j}{\partial {\mathbf n}} ds =
2\pi \int_{\Gamma_p}\frac{\partial G(z,z_j)}{\partial {\mathbf n}} ds -
\int_{\Gamma_p}\frac{\partial\log\abs{z-z_j}}{\partial {\mathbf n}} ds = $$
$$= 2\pi \int_{\Gamma_p}\frac{\partial G(z,z_j)}{\partial {\mathbf n}} ds = - 2\pi h_p(z_j)$$
for every $p=1,\dots,m$ and $j=1,\dots,M$.

This calculation and (\ref{conjugate}) give
\begin{equation}\label{surjective}
\frac{1}{2\pi}\int_{\Gamma_p}\frac{\partial u_R}{\partial {\mathbf n}} ds = 
\sum_{j=1}^M h_p(z_j) + k_p
\end{equation}
for every $p=1,\dots,m$. 
Since $R$ is an arbitrary positive $C^{k+1}$ $(k\ge 3)$ function, $u_R$ is an arbitrary real $C^{k+1}$ function
on $\bdr D$ and the left hand sides of (\ref{surjective}) are arbitrary real numbers.
Indeed, let $\lambda_1,\dots,\lambda_m$ be real constants such that
$$\sum_{p=1}^{m} \lambda_p \int_{\Gamma_p}\frac{\partial u}{\partial {\mathbf n}} ds = 
\int_{\bdr D} \frac{\partial u}{\partial {\mathbf n}} (\sum_{p=1}^m \lambda_p h_p) ds = 0$$
for all $C^{k+1}(D)$ functions $u$ which are harmonic on $D$. Then we get by the Green's identity
that
$$\int_{\bdr D} u (\sum_{p=1}^m \lambda_p \frac{\partial h_p}{\partial {\mathbf n}}) ds = 0$$
for every such $u$ (here ${\mathbf n}$ is the outer unit normal to $\bdr D$). Therefore
$$\sum_{p=1}^m \lambda_p \frac{\partial h_p}{\partial {\mathbf n}} = 0$$
on $\bdr D$ and hence the harmonic function $h = \sum_{p=1}^m \lambda_p h_p$ is constant by the
Hopf's lemma. Since $h=0$ on $\Gamma_0$ we get that $\lambda_1 = \lambda_2 = \dots = \lambda_m = 0$
and the left hand sides of (\ref{surjective}) are arbitrary real numbers.
Therefore we get surjectivity (in the case $M<m$, points $z_{M+1},\dots,z_m$ are chosen from $\Gamma_0$).
\end{proof}

\begin{example}\label{example}
Complete computations can be carried over for an annulus $A(q,1) = \D(0,1)\setminus\overline{\D(0,q)}$.
In this case 
$$h_1(z) = \frac{\log(\abs{z})}{\log(q)}$$
and $h_1$ obviously takes all the values from the interval $\lbrack 0,1\rbrack$. 
Then the equation (\ref{surjective}) is
\begin{equation}\label{lhs}
\frac{1}{2\pi}\int_{\Gamma_1}\frac{\partial u_R}{\partial {\mathbf n}} ds = 
\sum_{j=1}^M h_1(z_j) + k_1
\end{equation}
and we can easily get a solution of the Riemann-Hilbert problem with either $1$
(the left hand side of (\ref{lhs}) is not an integer) or no zeros (the left hand side of
(\ref{lhs}) is an integer).
Observe also, that there is always a solution of the Riemann-Hilbert problem
with $2$ zeros. Finally we observe that local 'pushing' of solutions
of the Riemann-Hilbert problem using the implicit function theorem in general 
has to lead to a situation
in which the index of the problem is such that no further direct
use of the implicit function theorem is possible.
\end{example}

\section{Appendix}
\label{appendix}

In the appendix we give proofs of some technical results needed in our 
arguments.

\begin{lemma} \label{rational}
Let $T$ be a finite disjoint union of $C^k$ $(k\ge 2)$ maximal real tori in $\bdr\D\times\C$. Then for
each point $(z_0,w_0)\in(\bdr\D\times\C)\setminus T$ there exists a polynomial $P$ such that 
$P(z_0,w_0)=0$ and that the intersection of its zero set $V=\zp (z,w)\in\C^2; P(z,w)=0\zz$ with
$\bdr\D\times\C$ is a smooth simple closed curve which does not intersect $T$.
In particular $T$ is a rationally convex subset of $\C^2$.
\end{lemma}
\begin{proof}
If $T$ is a maximal real torus in $\bdr\D\times\C$ such that each fiber $T\cap (\zp z\zz\times\C)$
has only one connected component, then for each point $(z_0,w_0)\in	(\bdr\D\times\C)\setminus T$
there exists a smooth closed curve $z\mapsto (z,w(z))$ (a mapping from $\bdr\D$ into $\bdr\D\times\C$)
such that $(z,w(z))\in(\bdr\D\times\C)\setminus T$ for every $z\in\bdr\D$ and $w(z_0)=w_0$. Approximating
the smooth function $w(z)$ uniformly on $\bdr\D$ by Fourier series we may assume 
$$w(z) = \sum_{j=-N}^{j=N} a_j z^j = \frac{p(z)}{z^N}$$
for some holomorphic polynomial $p(z)$.
We define $P(z,w) = z^N w - p(z)$.
Similar argument works in the case $T=\cup_{j=1}^n T_j$ is a finite disjoint union of maximal real tori $T_j$ in
$\bdr\D\times\C$ such that each component $T_j$ has connected fibers.

Let now $T$ be a maximal real torus in $\bdr\D\times\C$ such that each fiber $T\cap (\zp z\zz\times\C)$
has $k$ connected components. Let $\Phi(z,w) = (z^k,w)$ be a proper polynomial map from 
$\overline{\D}\times\C$ into itself.
Then the preimage $\Phi^{-1}(T)$ is a disjoint union of $k$ maximal real tori in $\bdr\D\times\C$ such that
each component of $\Phi^{-1}(T)$ has connected fibers. 

Let $(z_0,w_0)\in (\bdr\D\times\C)\setminus T$ and let $(z_1,w_1)\in \Phi^{-1}(z_0,w_0)$.
Then $(z_1,w_1)\in (\bdr\D\times\C)\setminus \Phi^{-1}(T)$ and we can find a polynomial
$Q(z,w) = z^N w - q(z)$ such that its zero set $V$ passes through $(z_1,w_1)$ and it intersects
$\bdr\D\times\C$ in a smooth simple closed curve which does not intersect $\Phi^{-1}(T)$.
The variety $\Phi(V)$ passes through the point $(z_0,w_0)$ and is the zero set of a polynomial $P(z,w)$.
Also, $\Phi(V)$ intersects $\bdr\D\times\C$ in a smooth simple closed curve which does not intersect $T$.
Similar argument works in the case $\Phi(z,w) = (z^{k\,n},w)$.

The general case can be treated similarly as the previous one. For each connected component
$T_0$ of $T$ let $k_0$ be the number of connected components in its fibers. We 'pull-back' the manifold $T$ 
with the map $\Phi(z,w)=(z^k,w)$ where $k$ is the product of the numbers $k_0$. The rest is as in the
previous paragraph.	
\end{proof}

From the work of Duval, \cite{duv}, we get the following lemma.

\begin{lemma}\label{symplectic}
Let $\zp\gamma_z\zz_{z\in\bdr\SI}$ be a $C^k$ $(k\ge 2)$ family of Jordan curves in $\C$ 
which all contain the point $0$ in their interior and let
$$T = \bigcup_{z\in\bdr\SI} (\zp z\zz\times\gamma_z)$$
be a disjoint union of totally real tori in $\bdr\SI\times\C$. 
Then there exists a $C^k$ strongly plurisubharmonic 
function $v$ on $\overline{\SI}\times\C$ such that
$T$ is Lagrangian for the symplectic form $\omega = i \partial\overline{\partial} v$.
In addition, we can arrange that the $\omega$-area of any disc $\zp z\zz\times\widehat{\gamma_z}$ is $1$.
\end{lemma}
\begin{remark}
One may always consider $\overline{\SI}$ as a subset of some larger open Riemann surface $\SI_1$. A function
$v$ is a $C^k$ strongly plurisubharmonic function on $\overline{\SI}\times\C$ if it is such on an open
neighbourhood of $\overline{\SI}\times\C$ in $\SI_1\times\C$.
\end{remark}
\begin{proof}
Let $f : \SI\rightarrow\D$ be a nonconstant holomorphic function, smooth up to the boundary, 
such that $\abs{f(z)} = 1$ for every $z\in\bdr\SI$, \cite{ahl}, and let
$g$ be a holomorphic function on $\SI$, smooth up to the boundary, such that the mapping
$z\mapsto (f(z),g(z))$ from $\overline{\SI}$ into $\overline{\D}\times\C$ is an immersion which 
is injective on $\bdr\SI$, \cite{cer-for}.

We observe the mapping $\Phi :\overline{\SI}\times\C\rightarrow\overline{\D}\times\C$ defined as
$\Phi(z,w) = (f(z), g(z) + w)$. We may assume that the mapping $\Phi$ is injective on $T$.
If not, replace $g$ with $R\,g$, for some $R>0$ large enough. 

The totally real manifold $\Phi(T)$ is a subset of $\bdr\D\times\C$ and it is rationally convex, 
Lemma \ref{rational}. Using the result of Duval, \cite{duv}, there exists a strongly plurisubharmonic 
function $u$ on $\C^2$ such that $\Phi(T)$ is Lagrangian for the symplectic 
form $\Omega = i\partial\overline{\partial} u$. In a neighbourhood of $\Phi(T)$ function $u$ can be 
defined as $u(z,w) = (\abs{z}^2-1)^2 + \rho(z,w)^2$, where $\rho$ is a defining function for $\Phi(T)$.
In addition, we may assume that outside some large ball
$u(z,w) = \lambda(\abs{z}^2 + \abs{w}^2)$ for some $\lambda>0$.

From the Stokes' theorem it follows that for any connected component 
$T_0$ of $\Phi(T)$ is the $\Omega$-area of any component of its fibers 
$T_0\cap (\zp z\zz\times\C)$, $z\in\bdr\D$, a constant which depends only on $T_0$.
We may assume that all these constants, which could be different for different components $T_0$, 
are smaller than $1$ (multiplay $u$ with a suitable positive constant).

Let $P(z,w)$ be a polynomial such that the intersection of its zero set with
$\bdr\D\times\C$ is a smooth curve contained in the bounded component of $(\bdr\D\times\C)\setminus T_0$
(Lemma \ref{rational}). Let $a>0$ be so small that the intersections of all level sets 
$P(z,w)=c$, $\abs{c}\le a$, with $\bdr\D\times\C$ are still contained in the bounded component of 
$(\bdr\D\times\C)\setminus T_0$.
Let $\varphi$ be a smooth nonnegative function on $\R$ such that its support lies in $\lbrack -a^2,a^2\rbrack$.
We observe the $(1,1)$ nonnegative closed form 
$$\nu= i\varphi(\abs{P}^2)\, \partial P\wedge\overline{\partial P}$$
on $\C^2$ whose support is contained in a neighbourhood of the zero set of $P$.
The pull-back of $\nu$ to $\zp z_0\zz\times\C$ is 
$$i\varphi(\abs{P(z_0,w)}^2)\, \abs{P_w(z_0,w)}^2 dw\wedge d\overline{w} .$$
Since $P$ is not identically equal to $0$ on any vertical line $\zp z\zz\times\C$, $z\in\bdr\D$,
the derivative $P_w(z_0,w)$ can be zero only at finite number of points $w$. Therefore we can
find a suitable function $\varphi$ such that the $(\Omega+\nu)$-area of any component of the fibers
of $T_0$ is equal to $1$. This procedure can be done independently for any connected component of $\Phi(T)$.

Since $\Omega+\nu$ is a positive $(1,1)$ closed form on $\C^2$, 
there exists a strongly plurisubharmonic function $u_1$ such that 
$\Omega + \nu = i\partial\overline{\partial} u_1$. 
We may assume that outside some ball $u_1$ equals a positive multiple of $\abs{z}^2 + \abs{w}^2$.
In addition, $\Phi(T)$ is a Lagrangian submanifold for $\Omega + \nu$ and the $(\Omega + \nu)$-area of any 
component of its fibers is $1$.

Let $\mathcal{C}=\zp z_1,\dots,z_k\zz$ be a finite set of points on $\SI$ where $df=0$. 
Then $v_1=u_1\circ\Phi$ is a strongly plurisubharmonic function on 
$(\overline{\SI}\setminus\mathcal{C})\times\C$ and $v_{1w\overline{w}} > 0$ on $\overline{\SI}\times\C$.
Let $\chi$ be a smooth nonnegative function on $\overline{\SI}$ whose
support is contained in some small neighbourhood of the set $\mathcal{C}$
and which equals $1$ in some smaller neighbourhood of $\mathcal{C}$.
Finally we define $v = v_1 + \varepsilon \chi \abs{g}^2$ and 
function $v$ has all the required properties if $\varepsilon$
is chosen small enough.	
\end{proof}

Now we give the necessary apriori estimates for approximate solutions of the Riemann-Hilbert
problem. We are working on $\D$ and with a $C^{k+1}$  $(k\ge 3)$ family of Jordan
curves $\zp \gamma_z \zz_{z\in\bD}$ in $\C$ such that the point $0$ is in the interior of each of them.
In what follows $C$ will denote a universal constant which depends only on the data. Also, we will use
the fact that the Hilbert transform is a bounded linear operator on the Sobolev spaces $W^{k,p}(\bD)$,
$1<p<\infty$, and H\"{o}lder spaces $C^{k,\alpha}(\bD)$.

\begin{lemma}\label{compositum}
Let $\varphi$ be a $C^2$ function on $\bD\times\C$ and let 
$\zp f_n\zz_{n=1}^{\infty}$ be a uniformly bounded sequence of $C^{1,\alpha}$
functions on $\bD$:
$$\norm{f_n}_{\infty} \le R<\infty.$$
Then there exists a constant $C$ which depends on the $C^2$ norm of $\varphi$ on 
$\bD\times\overline{\D(0,R)}$ such that
$$\norm{\varphi(z,f_n(z))}_{\alpha}\le C(\norm{f_n}_{\alpha} + 1) ,$$
and
$$\norm{\varphi(z,f_n(z))}_{1,\alpha}\le C(\norm{f_n}^{2}_{1,\alpha} + 1) .$$
\end{lemma}

\begin{remark}
For the first inequality it is enough to assume that $\varphi$ is a $C^1$ function.
\end{remark}

\begin{proof}
We denote $\varphi(\theta,w) := \varphi(e^{i\theta},w)$ and $f_n(\theta) := f_n(e^{i\theta})$.
Let $M$ be the $C^2$ norm of $\varphi$ on $\bD\times\overline{\D(0,R)}$. 
Then
$$\abs{\varphi(\theta,f_n(\theta))-\varphi(\nu,f_n(\nu))}\le$$
$$\le\abs{\varphi(\theta,f_n(\theta))-\varphi(\nu,f_n(\theta))} + 
\abs{\varphi(\nu,f_n(\theta))-\varphi(\nu,f_n(\nu))}\le $$
$$\le M(\abs{\theta -\nu} + \abs{f_n(\theta) - f_n(\nu)})$$
by Lagrange theorem. So 
$$\norm{\varphi(\theta,f_n(\theta))}_{\alpha} \le C(\norm{f_n}_{\alpha} + 1) .$$
Using this inequality for partial derivatives $\partial_{\theta}\rho$ and $D_w\rho$ we get
$$\norm{\varphi(\theta,f_n(\theta))}_{1,\alpha} = \norm{\varphi(\theta,f_n(\theta))}_{\infty}
+ \norm{\frac{\partial}{\partial\theta}(\varphi(\theta,f_n(\theta)))}_{\alpha}\le$$
$$\le M + \norm{(\partial_{\theta}\varphi)(\theta,f_n(\theta))}_{\alpha} +
\norm{(D_w\varphi)(\theta,f_n(\theta))}_{\alpha}\,\norm{\frac{\partial}{\partial\theta}f_n(\theta)}_{\alpha}\le$$
$$\le M + C(\norm{f_n}_{\alpha} + 1) + C(\norm{f_n}_{\alpha} + 1)\,\norm{f_n}_{1,\alpha}\le 
C (\norm{f_n}^2_{1,\alpha} + 1) .$$ 
\end{proof}

\begin{remark}
In general one can get an estimate 
$$\norm{\varphi(\theta,f_n(\theta))}_{k,\alpha}\le C(\norm{f_n}^{k+1}_{k,\alpha} + 1) .$$
\end{remark}
In the next lemma we closely analyze the proof of Forstneri\v{c} in \cite{for1}.
\begin{lemma}\label{estimates} There exist a finite positive constant $C$ and a positive integer $l$ 
such that for each $n\in\N$ and $0<q<1$ there is a solution $f_n$ of the Riemann-Hilbert problem on the disc 
with the winding number $n$ and such that 
$$\norm{f_n}_{A^{1,\alpha}(\D(0,1))} \le C n^l
\ \ \ \ {\rm and}\ \ \ \ \norm{f_n}_{A^{1,\alpha}(\D(0,q))}\le C n^l q^{n-2}.$$
\end{lemma}
\begin{proof} Let $\rho(e^{i\theta},w)=\rho(\theta,w)$ be a $C^{k+1}$ defining function for $\zp \gamma_z\zz_{z\in\bD}$.
We define 
$$\rho_n(\theta,w) = \rho(\theta, e^{in\theta} w)$$
and
$$\gamma_z^n = \zp w\in\C; \rho_n(\theta, w) = 0\zz = \frac{1}{z^n}\gamma_z$$
where $z = e^{i\theta}$.
Then there exist universal constants $0<r<R<\infty$ such that 
$\D(0,r)\subseteq\gamma_z^n\subseteq \D(0,R)$.

Let $\nu(\theta,w) = (\overline{\partial}_w\rho)(\theta,w)$ be the normal to
$\gamma_z$ $(z=e^{i\theta})$ at the point $w$ and $\eta(\theta,w) = w\overline{\nu(\theta,w)}$.
Similarly we define $\nu_n(\theta,w)$ and $\eta_n(\theta,w)$ for $\rho_n$.
Then 
$$\nu_n(\theta,w) = \nu(\theta, e^{in\theta} w) e^{-in\theta}\ \ \ \ {\rm and}\ \ \ \ 
\eta_n(\theta,w) = \eta(\theta,e^{in\theta} w) .$$
Geometric assumptions on $\zp \gamma_z \zz_{z\in\bD}$ imply that there exist well defined real $C^{k+1}$ functions
$a$ and $b$ on $T = \cup_{z\in\bD} \zp z\zz\times\gamma_z$ such that
$$\eta(\theta,w) = e^{a(\theta,w) + i b(\theta,w)} .$$
Hence 
$$\eta_n(\theta,w) =  e^{a_n(\theta,w) + i b_n(\theta, w)},$$
where $a_n(\theta,w) = a(\theta,e^{in\theta} w)$ and $b_n(\theta,w) = b(\theta,e^{in\theta}w)$.

Let $h_n = e^{g_n}$ be a solution of the Riemann-Hilbert problem for $T_n$ with no zeros on $\D$. 
We know, \cite{chi}, that $h_n$ is of class $\C^{k,\alpha}$ on $\overline{\D}$ for each $0<\alpha<1$.
We have the following immediate apriori uniform estimates:
$$\norm{h_n}_{\D}\le R ,\ \ \ \ \norm{a_n}_{T_n}\le \norm{a}_{T} , \ \ \ \  \norm{b_n}_{T_n}\le \norm{b}_{T} .$$
Also, 
\begin{equation}\label{partialthetarho}
(\partial_{\theta}\rho_n)(\theta,w) = 
(\partial_{\theta}\rho)(\theta, e^{in\theta}w) + 2{\rm Re}(\overline{\nu_n(\theta,w)} i n w)
\end{equation}
and so
$$\norm{(\partial_\theta \rho_n)(\theta,h_n(\theta))}_{\bD} 
\le \norm{\rho}_{C^1(T)} + 2 n R \norm{\rho}_{C^1(T)}\le C n$$
for some universal constant $C$. Here we used the notation $h_n(\theta) = h_n(e^{i\theta})$.

We differentiate the equation
$$\rho_n(\theta, h_n(\theta)) = 0$$
with respect to $\theta$ to get
$$(\partial_{\theta}\rho_n)(\theta, h_n(\theta)) + 
2{\rm Re}(\overline{\nu_n(\theta, h_n(\theta))}\, \frac{\partial h_n}{\partial\theta}(\theta)) = 0$$
and so 
\begin{equation}\label{basiceq}
{\rm Re}(\frac{\partial g_n}{\partial\theta}(\theta) e^{-\widetilde{b}_n(\theta) + i b_n(\theta, h_n(\theta))})=
- e^{\widetilde{b}_n(\theta)} e^{-a_n(\theta, h_n(\theta))} (\partial_{\theta}\rho_n)(\theta, h_n(\theta)),
\end{equation}
where $\widetilde{b}_n(\theta)$ is the Hilbert transform of $b_n(\theta,h_n(\theta))$ with the value 
$0$ at $0$. 

Let $2<p<\infty$. 
As it is proved in \cite{for1} on pages 881-882 we can write
$b = {\rm Re}(q) + b^{\prime}$ so that $p\,\norm{b^{\prime}}_{T} < \frac{\pi}{2}$
and $q$ is a function of the form
$$q(\theta,z) = \sum_{0\le k\le k_0}\sum_{-j_0\le j\le j_0} c_{jk} e^{ij\theta} z^k.$$
Then by a classical estimate \cite[p.114]{gar1} we have
$$\norm{e^{\pm\widetilde{b^{\prime}}}}_p \le C(p)$$
for some universal constant $C(p)$ which depends only on $p$.
Also, the operator $R_q$ which assigns to each $f\in C(\bD)$ the harmonic conjugate of
${\rm Re}(q(\theta, f(e^{i\theta})))$ is a bounded nonlinear operator from
$A(\bD)$ into $C(\bD)$. 

Thus we have the estimate $p\,\norm{b_n^{\prime}}_{T_n} < \frac{\pi}{2}$ for
$b_n^{\prime}(\theta,w) = b^{\prime}(\theta, e^{in\theta} w)$, and hence
$$\norm{e^{\pm\widetilde{b_n^{\prime}}}}_p \le C(p) .$$
We still have to estimate $R_q(z^n h_n(z))$. Since the family of holomorphic functions
$z^n h_n(z)$ is uniformly bounded and since the operator $R_q$ is a
bounded operator from $A(\bD)$ into $C(\bD)$ we get a uniform bound on the harmonic 
conjugates of ${\rm Re}(q(\theta, e^{in\theta} h_n(e^{i\theta})))$.
Therefore
$$\norm{e^{\pm \widetilde{b}_n}}_p = \norm{e^{\pm R_q(z^n h_n(z))}\, e^{\pm\widetilde{b_n^{\prime}}}}_p\le C(p) .$$
Now we can estimate the righthand side of (\ref{basiceq}) in $L^p$ norm
$$\norm{e^{\widetilde{b}_n(\theta)}\, e^{-a_n(\theta, h_n(\theta))}\, (\partial_{\theta}\rho_n)(\theta, h_n(\theta))}_p
\le C(p)\, n $$
and, since the Hilbert transform is bounded in $L^p$ spaces $(1<p<\infty)$, we get 
$$\norm{\frac{\partial g_n}{\partial\theta}(\theta)\, e^{-\widetilde{b}_n(\theta) + i b_n(\theta, h_n(\theta))}}_p
\le C(p)\, n .$$
Also, since 
$$\frac{\partial h_n}{\partial\theta} = h_n\, \frac{\partial g_n}{\partial\theta} = 
h_n\, e^{\widetilde{b_n} - i b_n}\, \frac{\partial g_n}{\partial\theta}\, e^{-\widetilde{b_n} + i b_n} $$
we get 
$$\norm{\frac{\partial h_n}{\partial\theta}}_{\frac{p}{2}}\le C(p)\, n$$
and so
$$\norm{h_n}_{1,\frac{p}{2}}\le C(p)\, n$$
for every $n$. Using the embedding of the space $W^{1,\frac{p}{2}}(\bD)$ into the space $C^{\alpha}(\bD)$
for $\alpha = 1 - \frac{2}{p}\in (0,1) $ we get
$$\norm{h_n}_{\alpha} \le C(\alpha)\, n.$$
Now we estimate the derivatives
$$\norm{\frac{\partial}{\partial\theta} (e^{\pm\widetilde{b_n}})}_{\frac{p}{2}} =
\norm{e^{\pm\widetilde{b_n}}\,\frac{\partial \widetilde{b_n}}{\partial\theta}}_{\frac{p}{2}}\le
\norm{e^{\pm\widetilde{b_n}}}_p\, \norm{\frac{\partial \widetilde{b_n}}{\partial\theta}}_p\le $$
$$\le C(p)\, 
\norm{\widetilde{b_n}}_{1,p}\le C(p)\, \norm{b_n}_{1,p} .$$
Since $\frac{\partial}{\partial\theta}(b_n(\theta, h_n(\theta)))$ equals
$$(\partial_{\theta} b)(\theta,e^{in\theta} h_n(\theta)) + 
2{\rm Re}((\partial_w b)(\theta, e^{in\theta} h_n(\theta))\,( e^{in\theta}\frac{\partial h_n}{\partial\theta}(\theta)
+ i n e^{in\theta} h_n(\theta)))$$
we get
$$\norm{b_n}_{1,p}\le C(p)\, n$$
and so $\norm{b_n}_{\alpha}\le C(\alpha)\, n$ and
$$\norm{e^{\pm\widetilde{b_n}}}_{1,\frac{p}{2}}\le C(p)\, n\ \ \ \ {\rm and}\ \ \ \ 
\norm{e^{\pm\widetilde{b_n}}}_{\alpha}\le C(\alpha)\, n\ \ \ \ (\alpha = 1 - \frac{2}{p}).$$

Now we can work by induction.
From (\ref{basiceq}) we get
$$\norm{\frac{\partial g_n}{\partial\theta}(\theta)\, e^{-\widetilde{b}_n(\theta) + 
i b_n(\theta, h_n(\theta))}}_{\alpha}
\le\norm{e^{-\widetilde{b_n}}}_{\alpha}\, \norm{e^{-a_n}}_{\alpha}\,
\norm{(\partial_{\theta}\rho_n)(\theta,h_n(\theta))}_{\alpha}$$
and hence, using (\ref{partialthetarho}) and Lemma \ref{compositum},
$$\norm{\frac{\partial g_n}{\partial\theta}(\theta)\, e^{-\widetilde{b}_n(\theta) + 
i b_n(\theta, h_n(\theta))}}_{\alpha}
\le C n\, \norm{e^{in\theta} h_n(\theta)}_{\alpha}\, (n \norm{e^{in\theta} h_n(\theta)}_{\alpha})\le C\, n^6.$$
Since
$$\norm{\frac{\partial h_n}{\partial\theta}}_{\alpha} \le \norm{h_n}_{\alpha}\,
\norm{e^{\widetilde{b_n} - i b_n}}_{\alpha}\,
\norm{\frac{\partial g_n}{\partial\theta}\, e^{-\widetilde{b_n} + i b_n}}_{\alpha},$$
we get 
$$\norm{h_n}_{1,\alpha}\le  C\, n^{9}.$$

Finally we define $f_n(z) = z^n h_n(z)$. Then 
$$\norm{f_n}_{A^{1,\alpha}(\D(0,1))} \le C\,n^{11}
\ \ \ \ {\rm and}\ \ \ \ \norm{f_n}_{A^{1,\alpha}(\D(0,q))}\le C\, n^{11}\, q^{n-2}$$
for any $n\in\N$ and $0<q<1$. 
\end{proof}

\begin{remark}
It can be shown by induction that for each fixed $k\in\N$ and $0<\alpha<1$ the norms of 
$\norm{h_n}_{k,\alpha}$ increase as some fixed power of $n$.
\end{remark}

\begin{lemma}\label{nonhomogeneous}
There exist constants $C>0$ and $l_1\in\N$ such that 
the linear operator $\widetilde{\Psi_n} : A^{1,\alpha}(\D)\rightarrow C^{1,\alpha}(\bdr\D)$
defined as
$$(\widetilde{\Psi_n} h)(z) = 2{\rm Re}((\partial_w\rho)(z,f_n(z)) h(z)) .$$
has a bounded right inverse $\widetilde{B_n} : C^{1,\alpha}(\bdr\D)\rightarrow A^{1,\alpha}(\D)$ with the norm
$\norm{\widetilde{B_n}}\le C\, n^{l_1}$ and such that 
$$\norm{\widetilde{B_n} g}_{A^{1,\alpha}(\D(0,q))}\le C n^{l_1} q^{n-2} \norm{g}_{1,\alpha}$$
for every $g\in C^{1,\alpha}(\bdr\D)$.
\end{lemma}

\begin{proof}
Let $g\in C^{1,\alpha}(\bdr\D)$. We define
$$\widetilde{B_n} g = \frac{1}{2}\,f_n e^{(-i)(b_n + i \widetilde{b_n})} 
( e^{-a_n-\widetilde{b_n}}\,g + i\, T(e^{-a_n-\widetilde{b_n}}\,g))$$
where $T : C^{1,\alpha}(\bdr\D)\rightarrow C^{1,\alpha}(\bdr\D)$ is the Hilbert transform.
From the estimates in the proof of the previous lemma and Lemma \ref{compositum} we have
$$\norm{e^{- a_n}}_{1,\alpha}\le C\,(\norm{f_n}^2_{1,\alpha}+1)\le C n^{2l},\ \ \ \  
\norm{e^{-i b_n}}_{1,\alpha}\le C\,(\norm{f_n}^2_{1,\alpha}+1)\le C n^{2l}$$
and
$$\norm{e^{\pm \widetilde{b_n}}}_{1,\alpha} \le \norm{e^{\pm\widetilde{b_n}}}_{\bD} + 
\norm{e^{\pm\widetilde{b_n}}\frac{\partial\widetilde{b_n}}{\partial\theta}}_{\alpha}
\le C n + \norm{e^{\pm\widetilde{b_n}}}_{\alpha}\norm{\widetilde{b_n}}_{1,\alpha}\le C n^{2l+1},$$
where the fact that the Hilbert transform is bounded on $C^{1,\alpha}(\bD)$ was used.
Then $\widetilde{\Psi_n}\widetilde{B_n} = I$ and
$$\norm{\widetilde{B_n} g}\le C\, \norm{f_n}_{1,\alpha}\,\norm{e^{-i b_n}}_{1,\alpha}\,
\norm{e^{\widetilde{b_n}}}_{1,\alpha}\,
\norm{e^{-a_n}}_{1,\alpha}\,\norm{e^{-\widetilde{b_n}}}_{1,\alpha}\,\norm{g}_{1,\alpha} .$$
Hence $\norm{\widetilde{B_n} g}\le C\, n^{9l+2}\, \norm{g}_{1,\alpha}$ and
similarly
$$\norm{\widetilde{B_n} g}_{A^{1,\alpha}(\D(0,q))}\le C n^{9l+2}\, q^{n-2}\, \norm{g}_{1,\alpha} .$$
\end{proof}

With the notation of the proof of Theorem \ref{main} we state and prove the following
lemma.
\begin{lemma}\label{approximation}	The estimate
$$\norm{\rho(z, h_n(z))}_{1,\alpha}\le C n^{3l} q^{n}$$
holds for some universal finite constant $C$ and every $n\in\N$.
\end{lemma}

\begin{proof}
Let $\varphi(t) =\rho(z, t h_n(z) + (1-t) f_n(z))$ be a mapping from
the unit interval $\lbrack 0,1\rbrack$ into $C^{1,\alpha}(\bD)$. Then
$$\rho(z, h_n(z)) = \rho(z, h_n(z)) - \rho(z, f_n(z))	=
\varphi(1)-\varphi(0) =$$
$$=\int_0^1 \varphi^{\prime}(t)\,dt = 
\int_0^1 2{\rm Re}((\partial_w\rho)(z, t h_n(z) + (1-t) f_n(z))(h_n(z)-f_n(z)))\,dt .$$
Hence
$$\norm{\rho(z, h_n(z))}_{1,\alpha}\le 
\int_0^1 2\norm{(\partial_w\rho)(z, t h_n(z) + (1-t) f_n(z))}_{1,\alpha}
\norm{h_n-f_n}_{1,\alpha}\,dt. $$
We observe that 
$$\norm{h_n-f_n}_{1,\alpha}\le C\, n^{l} q^{n}$$ and from Lemma \ref{compositum}
$$\norm{(\partial_w\rho)(z, t h_n(z) + (1-t) f_n(z))}_{1,\alpha}\le
C\, (\norm{t h_n + (1-t) f_n}_{1,\alpha}^2 +1)\le$$
$$\le  C\, (\norm{f_n}_{1,\alpha}^2+\norm{h_n-f_n}_{1,\alpha}^2+1) \le C n^{2l} .$$
Hence
$$\norm{\rho(z, h_n(z))}_{1,\alpha}\le C\, n^{3l} q^{n} .$$
\end{proof}
{\bf Acknowledgement.}
I wish to thank Manolo Flores for stimulating discussions.

\end{document}